\documentclass{article}
\usepackage{amssymb,amsmath}

\newtheorem{theorem}{Theorem}
\newtheorem{corollary}{Corollary}

\begin{document}
\title{Structure of the time projection\\
for stopping times in von Neumann algebras}
\author{Andrzej \L uczak\thanks{Work supported by KBN grant 2 PO3A 04410}\\
 Faculty of Mathematics\\
 \L\'od\'z University\\
 ul. Banacha 22\\
 90-238 \L od\'z, Poland\\
\texttt{anluczak@math.uni.lodz.pl}}
\date{}
\maketitle

\begin{abstract}
We give an explicit formula for the time projection in an
arbitrary von Neumann algebra from which all its basic properties
can be easily derived. The analysis of the situation when this
time projection is a conditional expectation is also performed.
\end{abstract}

\section{Introduction}
The aim of these notes is to investigate
some properties of the time projection for a stopping time in a
von Neumann algebra. This is done solely by using an explicit
formula for the projection, without any reference to stochastic
integration. In particular, we obtain simple conditions for
stopping a noncommutative martingale. The problem of when the time
projection can be treated as a conditional expectation is also
addressed. Its solution, known in the case of the Clifford
probability gauge space, is thus generalised to a fairly general
context.

\section{Preliminaries and notation}Throughout the paper $\mathcal A$ will denote a
von Neumann algebra acting in a Hilbert space $\mathcal H$ with a
cyclic and separating vector $\Omega$. $\omega$ will stand for a
(normal faithful) vector state on $\mathcal A$ induced by
$\Omega$. Let $(\mathcal A_t : t\geq 0)$ be a filtration of
$\mathcal A$, i.e. an increasing net of von Neumann subalgebras of
$\mathcal A$ such that $\mathcal A=\mathcal
A_{\infty}:=(\cup_{t\geq 0}\mathcal A_t)''$. We assume that there
are normal conditional expectations $\mathbb E_{t}, \ t\geq 0$,
from $\mathcal A$ onto $\mathcal A_t$ leaving $\omega$ invariant.
It follows easily (cf.\cite[Proposition 1.2]{BT1}) that if we
define
\[
P_t (x\Omega )=(\mathbb E_{t}x)\Omega\ , \ x\in \mathcal A,
\]
then $P_t$ is a projection from $\mathcal H$ onto $\mathcal H_t =
\overline {\mathcal A_{t} \Omega}$, consequently, $P_t \in
\mathcal A_{t}'$; we have also $\mathbb E_t \mathbb E_s = \mathbb
E_s \mathbb E_t =\mathbb E_{s\land t}$. In what follows we shall
be concerned with the ``time parameter'' $t$ belonging either to
the interval $[0,+\infty)$ or to the interval $[0,u]$, where
$0<u\leq +\infty $. Accordingly, we adopt the following
definition. A (quantum, noncommutative) stopping time $\tau$ is an
increasing net $(q_t)$, $t\in [0,+\infty )$ or $[0,+\infty]$ of
projections such that $q_t \in \mathcal A_t ,\ q_0=0$, and
$\bigvee_{t\geq 0}q_t=\boldsymbol{1}$ in the case $t\in
[0,+\infty)$ or $q_{\infty}=\boldsymbol{1}$ in the case $t\in
[0,+\infty]$. The definition above is a proper generalisation of
the notion of the classical (commutative) stopping time (cf.
\cite{BL,BT1,BT2,BW1} for more information). A fairly general
theory of stopping a noncommutative process has so far been
achieved only for martingales. Let us briefly recall its main
points here.

A martingale in $\mathcal H$ is a process $(\xi(t)\colon t\geq 0)$
such that $\xi (t) \in \mathcal H_t$ and for each $s,t\geq 0,\
s\leq t$,
\[
P_{s}\xi (t)=\xi (s).
\]
If we allow $t\in [0,+\infty]$ then it follows that there is $\xi
(=\xi(\infty))$ such that $\xi (t)=P_t \xi$; such martingales are
called closed, and it is not difficult to see that the following
conditions are equivalent: (i)\ $(\xi (t))$ is closed (ii)\
$\sup_t \|\xi (t) \| <+\infty $ (iii)\  there exists
$\lim_{t\to\infty} \xi(t)$ (cf.\cite[Proposition 1.1]{BW1}).

Now stopping $(\xi (t))$ consists in the following procedure. For
interval $[0,u]\ (u=+\infty \ \text{if}\  (\xi (t))\ \text{is
closed})$ we consider its partition $\theta =\{ 0=t_0<t_1< \dots
<t_n=u\}$, and form the sum
\begin{equation}\label{e1}
\xi_{\tau (\theta)}=\sum\limits_{i=1}^{n} (q_{t_{i}}-q_{t_{i-1}})
\xi (t_i).
\end{equation}
Taking the limit of the net $\{\xi_{\tau (\theta)}\colon \theta
\text{ -- partition}\}$ as $\theta$ refines, gives us the stopped
element $\xi_{\tau}(u)$, which is all we need if $u=+\infty$;
however, if $u<+\infty$ it seems reasonable to define $\xi_{\tau}$
as $\lim_{u\to\infty} \xi_{\tau}(u)$.

The existence of the two limits above is by no means obvious. It
turns out that while the limit in \eqref{e1} does exist it need
not be so with the other one, and thus we are guaranteed only of
the possibility of stopping a closed martingale. To analyse
$\xi_{\tau(\theta)}$ observe that the martingale property yields
\[
\xi (t_i)=P_{t_i} \xi (u), \ i=1,\dots , n,
\]
and hence
\begin{equation}\label{e2}
\xi_{\tau(\theta)}=\sum_{i=1}^{n} (q_{t_{i}}-q_{t_{i-1}}) P_{t_i}
\ \xi (u).
\end{equation}
Put
\[
M_{\tau(\theta)}(u)=\sum_{i=1}^{n}(q_{t_{i}}-q_{t_{i-1}})P_{t_i}.
\]
Then $M_{\tau(\theta)}(u)$ is a projection in $\mathcal H$ (recall
that $P_{t_i} \in \mathcal A_{t_i}', \ q_{t_{i-1}}, q_{t_i}\in
\mathcal A_{t_i}$). It is easily seen that the net
$\{M_{\tau(\theta)}(u)\colon \theta\text{ -- partition}\}$
decreases, so there exists $\lim_{\theta} M_{\tau(\theta)}(u)$
which we denote by $M_{\tau}(u)$ and call the time projection; it
is also clear that
\[
M_{\tau}(u)=\bigwedge_{\theta}\ M_{\tau(\theta)}(u).
\]
Accordingly, we have by \eqref{e2}
\[
\xi_{\tau}(u)=\lim_{\theta}\
\xi_{\tau(\theta)}=\lim_{\theta} M_{\tau(\theta)}(u)\xi
(u)=M_{\tau}(u)\xi(u).
\]
If $u=+\infty$ we shall write $M_{\tau}$ instead of
$M_{\tau}(\infty)$; note that this is the case considered in
\cite{BL,BT2,BW1} and mainly in \cite{BT1}. However, in \cite{BT1}
a more general setting that we have defined above is also taken
into account.

As a final remark let us observe that the definition of the time
projection as well as the results of the next section could be
obtained for Haagerup's $L^2(\mathcal A,\omega)$-space and the
algebra $\mathcal A$ acting on it by left multiplication,
especially in view of a spatial isomorphism between the
representations $(\mathcal A,\mathcal H,\Omega)$ and $(\mathcal A,
L^2 (\mathcal A, \omega ),h_{\omega}^{1/2})$ where
$h_{\omega}^{1/2}$ is a cyclic and separating vector in
$L^2(\mathcal A,\omega)$. The reasons for which we have adopted a
more traditional approach lie in Section 3. There we want to treat
the time projection, which is a projection in a Hilbert space, as
a projection in the algebra $\mathcal A$, and passing from one to
another is much more straightforward in our original setup where
we have a natural embedding of $\mathcal A$ into $\mathcal H$
given by $\mathcal A\ni x\mapsto x\Omega\in \mathcal H$.

\section{Representation of the time projection} In this section we analyse
various properties of the time projection by means of an explicit
formula expressing it in terms of the $P_{t}$ and $q_{t}$.

\begin{theorem}\label{t1}Let $u\in (0, +\infty]$. Then
\begin{equation}\label{e3}
M_{\tau}(u)=\ \bigwedge_{t \leq u}\ (q_{u} - q_{t} P_{t}^{\perp}).
\end{equation}
\end{theorem}

\emph{Proof.} Take the partition $\theta_0 = \{0=t_0 < t_1 =u\}$.
We have
\[
M_{\tau}(u) \leq M_{\tau (\theta_{0})}=(q_{t_{1}} -
q_{t_{0}})P_{t_{1}}=q_u P_u \leq q_u .
\]

Let $\xi \in \mathcal H$, and assume that $M_{\tau}(u) \xi =q_u \xi$. For an
arbitrary $t\in [0, u]$ we have
\[
M_{\tau}(u) \leq (q_{t} - q_{0})P_t + (q_u - q_t)P_u \leq q_u ,
\]
giving the equality
\[
(q_t - q_0 )P_t \xi + (q_u - q_t )P_u \xi =q_u \xi .
\]
Applying $q_t$ to both sides yields
\begin{equation}\label{e4}
q_t P_t \xi \ =\ q_t \xi.
\end{equation}
Conversely, if for each $t\in [0, u]$ equality \eqref{e4} holds,
then for any $s\leq t$ we have, applying $q_s$ to both sides of
\eqref{e4},
\[
q_s P_t \xi = q_s \xi ,
\]
and for any partition $\theta =\{0=t_0 < t_1 < \dots < t_n =u \}$
\[
M_{\tau(\theta)}(u) \xi
=\sum_{i=1}^{n}(q_{t_{i}}-q_{t_{i-1}})P_{t_i} \xi =
\sum_{i=1}^{n}(q_{t_{i}}\xi -q_{t_{i-1}}\xi )=q_{t_n}\xi -
q_{t_0}\xi = q_u \xi ,
\]
hence
\[
M_{\tau}(u) \xi =\lim_{\theta} M_{\tau (\theta)} (u) \xi = q_u
\xi.
\]
We have thus obtained equivalence of the following conditions:
\begin{enumerate}
\item[(i)] $M_{\tau}(u)\xi = q_u \xi$
\item[(ii)] for each \ $t\in [0, u]\ \ q_t P_t \xi =q_t \xi$,
\end{enumerate}
or put in another way
\begin{enumerate}
\item[(i')] $[q_u - M_{\tau}(u)]\xi =0$
\item[(ii')] for each \ $t\in [0, u]\ \ q_t P_{t}^{\perp}\xi =0$.
\end{enumerate}

But condition (ii') is equivalent to the equality
\[
\bigl( \bigvee_{t\leq u}q_t P_{t}^{\perp}\bigr)\xi =0 ,
\]
which means that the projections $q_u - M_{\tau}(u)$ and
$\bigvee_{t\leq u} q_t P_{t}^{\perp}$ have the same null spaces,
so they must be equal:
\[
q_u - M_{\tau}(u) =\bigvee\limits_{t\leq u} q_t P_{t}^{\perp} .
\]
Consequently,
\[
M_{\tau}(u) =q_u - \bigvee_{t\leq u} q_t P_{t}^{\perp}=
\bigwedge_{t \leq u}\ (q_{u} - q_{t} P_{t}^{\perp}).
\]
\begin{corollary} If $u=+\infty$ then
\begin{equation}\label{e5}
M_{\tau} = \bigwedge_{t \geq 0}(q_{t}^{\perp} + q_t P_t).
\end{equation}
\end{corollary}

Indeed, we then have
\[
q_{\infty} - q_t P_{t}^{\perp} = \boldsymbol{1} - q_t P_{t}^{\perp}=q_{t}^{\perp} +
q_t P_t ,
\]
and for $t=+\infty$
\[
q_{\infty}^{\perp}+ q_{\infty}P_{\infty}=\boldsymbol{1},
\]
giving
\[
M_{\tau}=M_{\tau}(\infty)=\bigwedge_{0\leq t\leq +\infty}
(q_{\infty} - q_t P_{t}^{\perp})=\bigwedge_{0\leq t\leq +\infty}
(q_{t}^{\perp} + q_t P_t)=\bigwedge_{0\leq t< +\infty}
(q_{t}^{\perp} + q_t P_t).
\]

\begin{theorem}\label{t2}Let the set $\{ M_{\tau} \xi (t): \ t\in
[0,+\infty)\}$ be norm-bounded. Then the martingale $(\xi(t))$ can be
stopped and
\[
\xi_{\tau}=\lim_{t\to\infty}M_{\tau}\xi (t).
\]
\end{theorem}

\emph{Proof.} Put
\begin{equation}\label{e6}
\eta (t)=M_{\tau}\xi (t).
\end{equation}
For each $s, t\in [0,+\infty)$ we have
\[
P_s (q_{t}^{\perp} + q_t P_t)=
 \begin{cases} P_s q_{t}^{\perp}+P_s q_t=P_s \
 &\text{for $s\leq t$}\\ P_s q_{t}^{\perp}+q_tP_t \ &\text{for $s>t$}
 \end{cases}
=(q_{t}^{\perp} + q_t P_t)P_s,
\]
and from \eqref{e5} we get
\[
P_s M_{\tau}=M_{\tau}P_s .
\]
If $s\leq t$, then
\[
P_s \eta (t)=P_s M_{\tau} \eta (t)=M_{\tau} P_s \eta (s) =M_{\tau}
\eta (s)= \eta (s),
\]
which shows that $(\eta (t))$ is a martingale, and since it is
norm-bounded, we have $\eta (t) \to \eta ,\ \text{as} \quad t\to
\infty$, for some $\eta\in \mathcal H$. From \eqref{e6} we have
\[
M_{\tau}(t) \eta (t)=M_{\tau}(t)M_{\tau} \xi (t)=M_{\tau} (t)\xi
(t)=\xi_{\tau} (t).
\]

Now
\begin{align*}
&\| M_{\tau}(t)\eta (t)-M_{\tau}\eta \|\leq\|
M_{\tau} (t)[\eta (t)-\eta] \| + \| [M_{\tau}
(t)-M_{\tau}]\eta \| \leq \\
&\leq\| \eta (t)-\eta\| +
\| [M_{\tau}(t)-M_{\tau}]\eta\| \ \to 0 ,
\end{align*}
since $\lim_{t\to\infty} M_{\tau}(t)=M_{\tau}$, consequently
\[
\xi_{\tau}(t)=M_{\tau}(t)\eta (t)\to M_{\tau}\eta .
\]
But $M_{\tau}\eta (t)=\eta (t)$, and thus $M_{\tau}\eta=\eta$,
giving
\[
\xi_{\tau}=\lim_{t\to\infty}\xi_{\tau}(t)=M_{\tau}\eta=\eta
=\lim_{t\to\infty}M_{\tau}\xi (t) .
\]

Observe that the result of the last theorem perfectly agrees with what we
have for a closed martingale where also
\[
 \xi_{\tau} =  M_{\tau}\xi = \lim_{t\to\infty}M_{\tau} \xi (t) .
\]

\section{Time projection as a conditional expectation} In this section we
consider a question when the time projection can be treated as a
conditional expectation. A problem of this type was analysed in \cite{BL} for
the Clifford probability gauge space and solved by using some properties
of the Clifford quantum stochastic integral. The solution we give here
works in the general context of an arbitrary von Neumann algebra;
moreover it is simple and does not employ any theory of stochastic
integration.

Let $\tau =(q_t : \ t\in [0,+\infty])$ be a stopping time, and let
$M_{\tau}$ be the time projection. $M_{\tau}$ can be treated as a
conditional expectation if
\[
M_{\tau}(x\Omega)=y\Omega ,
\]
and the map $\mathbb E_{\tau}\colon x\mapsto y$ is a conditional
expectation. We then have
\[
(\mathbb E_{\tau}x)\Omega = M_{\tau}(x\Omega).
\]
Put
\[
\mathcal B_{\tau}=\{ x\in \mathcal A \colon \text{for each}\ t\geq 0\quad xq_t =q_t x\}
=\mathcal A \cap \{q_t : t\in [0,+\infty]\}' .
\]
For any partition $\theta =\{0=t_0<t_1< \dots <t_n =+\infty \}$
let
\[
\mathcal A_{\tau (\theta)}=\{x\in \mathcal A \colon xq_{t_i}=q_{t_i}x \in \mathcal A_{t_i},
i=0, 1, \dots , n \},
\]
and let
\[
\mathcal A_{\tau} = \bigcap_{\theta} \mathcal A_{\tau (\theta)}=\{
x\in \mathcal A \colon \text{for each}\  t\geq 0\ xq_t=q_tx \in
\mathcal A_t \}.
\]

\begin{theorem}\label{t3} $M_{\tau} | \mathcal B$ is a normal faithful
conditional expectation onto $\mathcal A_{\tau}$ leaving $\omega$ invariant.
\end{theorem}

\emph{Proof.} For a partition $\theta =\{0=t_0<t_1< \dots <t_n
=+\infty \}$ define on $\mathcal B_{\tau}$ the map $\mathbb
E_{\tau_{(\theta )}}$ by
\begin{align*}
&\mathbb E_{\tau_{(\theta )}} x=\sum_{i=1}^{n}(q_{t_i} -
q_{t_{i-1}}) \mathbb
E_{t_i}x= \sum_{i=1}^{n} \mathbb E_{t_i} ((q_{t_i} - q_{t_{i-1}})x) = \\
&=\sum_{i=1}^{n} (\mathbb E_{t_i}x) (q_{t_i} - q_{t_{i-1}}),\
x\in \mathcal B_{\tau}.
\end{align*}
For each $t\in [0, +\infty]$ we have $t_{j-i}\leq t < t_j$ with
some $j$, so
\begin{align*}
&q_t \mathbb E_{\tau (\theta)} x = q_t \sum_{i=1}^{j-1}
(q_{t_i}-q_{t_{i-1}}) \mathbb E_{t_i} x + q_t
(q_{t_j}-q_{t_{j-1}})\mathbb E_{t_j} x + \\
&+q_t \sum_{i=j+1}^{n} (q_{t_i}-q_{t_{i-1}}) \mathbb E_{t_i} x =
\sum_{i=1}^{j-1} (q_{t_i}-q_{t_{i-1}})\mathbb E_{t_i} x +
(q_t-q_{t_{j-1}}) \mathbb E_{t_j} x,
\end{align*}
and
\begin{align*}
&(\mathbb E_{\tau (\theta)}x)q_t = \sum_{i=1}^{j-1} (\mathbb
E_{t_i} x) (q_{t_i}-q_{t_{i-1}}) q_t + (\mathbb
E_{t_j}x)(q_{t_j}-q_{t_{j-1}})q_t +\\
& + \sum_{i=j+1}^{n}(\mathbb
E_{t_i} x)(q_{t_i}-q_{t_{i-1}})q_t =
\sum_{i=1}^{j-1}(\mathbb E_{t_i} x)(q_{t_i}-q_{t_{i-1}}) +
(\mathbb E_{t_j} x) (q_t-q_{t_{j-1}}).
\end{align*}
But for $x\in \mathcal B_{\tau}$
\[
(q_{t_i}-q_{t_{i-1}})\mathbb E_{t_i} x = \mathbb E_{t_i}
((q_{t_i}-q_{t_{i-1}})x) = (\mathbb E_{t_i}
x)(q_{t_i}-q_{t_{i-1}}),
\]
and
\[
(q_t-q_{t_{j-1}})\mathbb E_{t_j} x = \mathbb E_{t_j}
x((q_t - q_{t_{j-1}})x) = (\mathbb E_{t_j} x)(q_t-q_{t_{j-1}}),
\]
which shows that
\[
q_t \ \mathbb E_{\tau (\theta)}x = (\mathbb E_{\tau (\theta)}x)q_t ,
\]
i.e. $\mathbb E_{\tau (\theta)}x \in \mathcal B_{\tau}$.
Furthermore, for each $j=0, 1, \dots , n$
\[
q_{t_j} \mathbb E_{\tau (\theta)}x =
\sum_{i=1}^{j}(q_{t_i}-q_{t_{i-1}}) \mathbb E_{t_i} x =
\sum_{i=1}^{j}(\mathbb E_{t_i} x)(q_{t_i}-q_{t_{i-1}})= (\mathbb
E_{\tau (\theta)}x)q_{t_j},
\]
showing that $\mathbb E_{\tau (\theta)}x \in \mathcal A_{\tau
(\theta)}$. For $x\in\mathcal A_{\tau (\theta)}$ we have $\mathbb
E_{t_i} x = x$, hence
\[
\mathbb E_{\tau (\theta)}x + \sum_{i=1}^{n}(q_{t_i} - q_{t_{i-1}})x=x
\]
which means that $\mathbb E_{\tau (\theta)}$ is a projection from
$\mathcal B_{\tau}$ onto $\mathcal B_{\tau}\cap \mathcal
A_{\tau_{(\theta)}}$. If $x\in \mathcal B_{\tau}^{+}$, then
$(q_{t_i} - q_{t_{i-1}})x =(q_{t_i} - q_{t_{i-1}})x(q_{t_i} -
q_{t_{i-1}})\geq 0$, so
\[
\mathbb E_{\tau (\theta)}x=\sum_{i=1}^{n}(q_{t_i} -
q_{t_{i-1}})\mathbb E_{t_i}x=\sum_{i=1}^{n} \mathbb
E_{t_i}((q_{t_i} - q_{t_{i-1}})x)\geq 0,
\]
thus $\mathbb E_{\tau (\theta)}$ is positive. Since $\mathbb
E_{\tau (\theta)}\boldsymbol{1}=\boldsymbol{1}$, we infer that $\|
\mathbb E_{\tau (\theta)}\| =1$, and by virtue of \cite[Theorem
9.1 p.116]{S}, $\mathbb E_{\tau (\theta)}$ is a conditional
expectation. We have
\[
(\mathbb E_{\tau (\theta)}x)\Omega = M_{\tau_{(\theta)}}(x\Omega), \ x\in \mathcal
B_{\tau}.
\]

Put
\[
x_{\theta}=\mathbb E_{\tau (\theta)}x.
\]
Then $\{ x_{\theta} \}$ is a bounded net of elements in $\mathcal
A$, and for each $x' \in \mathcal A'$
\[
x_{\theta} (x'\Omega)=x'(x_{\theta}\Omega)=x' (\mathbb
E_{\tau (\theta)}x)\Omega =x' M_{\tau
(\theta)}(x\Omega)\to x' M_{\tau}(x\Omega).
\]
Thus the net $\{ x_{\theta}\}$ converges on the dense subspace
$\mathcal A' \Omega$ of $\mathcal H$, and since
\newline $\| x_{\theta} \| \leq \| x\| $, it follows that $\{ x_{\theta}\}$ converges
in the strong operator topology on $\mathcal A$, consequently, there is $y\in
\mathcal A$ such that $ x_{\theta}\to y$ strongly.

Let
\[
\mathbb E_{\tau}x=y = \lim_{\theta}x_{\theta}=\lim_{\theta}\mathbb
E_{\tau (\theta)}x, \ x\in\mathcal B_{\tau}.
\]
Clearly, $\mathbb E_{\tau}$ is a linear positive map on $\mathcal
B_{\tau}$, such that
\[
(\mathbb E_{\tau}x)\Omega = M_{\tau}(x\Omega), \ x\in \mathcal B_{\tau}.
\]
Since $\mathbb E_{\tau}\boldsymbol{1} =\boldsymbol{1}$, we have $\|
\mathbb E_{\tau}\| = 1$. For any partition $\theta$ and $x\in
\mathcal B_{\tau}$,
\[
\mathbb E_{\tau (\theta)}(\mathbb E_{\tau}x)\Omega = M_{\tau (\theta)}((\mathbb
E_{\tau}x)\Omega)=M_{\tau(\theta)}M_{\tau}(x\Omega)=M_{\tau}(x\Omega)=(\mathbb
E_{\tau}x)\Omega ,
\]
showing that $\mathbb E_{\tau (\theta)}\mathbb E_{\tau}=\mathbb
E_{\tau}$, since $\Omega$ is separating.
\newline
Accordingly, $\mathbb E_{\tau}x\in \mathcal A_{\tau (\theta)}$ for each
$\theta$, and it follows that $\mathbb E_{\tau}x\in \bigcap_{\theta}\mathcal
A_{\tau (\theta)}=\mathcal A_{\tau}$.
Furthermore, if $x\in \mathcal A_{\tau}$, then $\mathbb E_{\tau (\theta)}x=x$ for
each $\theta$, so
\[
\mathbb E_{\tau}x=\lim_{\theta}\mathbb E_{\tau (\theta)}x=x,
\]
which means that $\mathbb E_{\tau}$ is a projection onto $\mathcal A_{\tau}$, and
thus a conditional expectation. From the equality
\[
\omega\circ\mathbb E_{\tau (\theta)}=\omega ,
\]
we obtain
\[
\omega\circ\mathbb E_{\tau}=\omega ,
\]
which, since $\mathbb E_{\tau}$ is positive, implies faithfulness and
normality of $\mathbb E_{\tau}.$

Let us observe that in an entirely analogous way we can obtain a
corresponding result for the time projection $M_{\tau}(u)$.
\newline
Indeed, putting
\[
\mathcal B_{\tau}(u)=\{ x\in \mathcal A : \text {for each}\
t\leq u \quad xq_{t}=q_t x \}
\]
\[
\mathcal A_{\tau}(u)=\{ x\in \mathcal A : \text {for each}\
t\leq u \quad xq_t=q_tx\in \mathcal A_t \} ,
\]
we get that $M_{\tau}(u)|\mathcal B_{\tau}(u)$ is a conditional expectation
onto $q_u \ \mathcal A_{\tau}(u)q_u$.

\end{document}